\documentclass[12pt]{article}
\usepackage{amssymb}
\usepackage{color}
\definecolor{c}{rgb}{0.9,0.3,0.1}
\definecolor{b}{rgb}{0.1,0.3,0.9}

\begin{document}
\setlength{\baselineskip}{24pt}
\newtheorem{thm}{Theorem}[section]
\newtheorem{example}[thm]{Example}
\newtheorem{cor}[thm]{Corollary}
\newtheorem{notation}[thm]{Notation}
\newtheorem{notations}[thm]{Notations}
\newtheorem{definition}[thm]{Definition}
\newtheorem{lemma}[thm]{Lemma}
\newtheorem{claim}[thm]{Claim}
\newtheorem{mthm}[thm]{Meta-Theorem}
\newtheorem{prop}[thm]{Proposition}
\newtheorem{rem}[thm]{Remark}
\newtheorem{conj}[thm]{Conjecture}
\newtheorem{rems}[thm]{Remarks}
\renewcommand {\theequation}{\thesection.\arabic{equation}}
\newcommand{\qed}{\hfill\rule{2mm}{3mm}\vspace{4mm}}

\def\al{{\alpha}}\def\be{{\beta}}\def\de{{\delta}}\def\De{{\Delta}}
\def\ep{{\epsilon}}
\def\ga{{\gamma}}\def\Ga{{\Gamma}}\def\la{{\lambda}}\def\om{{\omega}}
\def\si{{\sigma}}\def\th{{\theta}}\def\ze{{\zeta}}

\def\Om{{\Omega}}
\def\<{\left<}\def\>{\right>}\def\({\left(}\def\){\right)}

\newfam\msbmfam\font\tenmsbm=msbm10\textfont
\msbmfam=\tenmsbm\font\sevenmsbm=msbm7
\scriptfont\msbmfam=\sevenmsbm\def\bb#1{{\fam\msbmfam #1}}
\def\BB{\bb B}\def\CC{\bb C}
\def\EE{\bb E\ }\def\HH{\bb H}\def\NN{\bb N}\def\PP{\bb P}
\def\RR{\bb R}\def\WW{\bb W}

\def\cB{{\cal B}}\def\cF{{\cal F}}\def\cG{{\cal G}}
\def\cH{{\cal H}}\def\cK{{\cal K}}\def\cL{{\cal L}}\def\cM{{\cal M}}\def\cP{{\cal P}}
\def \cS{{\cal S}}\def\cY{{\cal Y}}\def\cZ{{\cal Z}}

\title{\large \bf Some properties for superprocess under a stochastic flow}
\medskip

\author{Kijung Lee$\,^a$, Carl Mueller$\,^{a,1}$
and Jie Xiong$\,^{b,c,2}$
 \\[0.5cm]
\small{$\;^a\,$Department of Mathematics, University of Rochester}\\
\small{Rochester, NY 14627}\\
 \small{$\;^b\,$Department of
Mathematics,
University of Tennessee,}\\
\small{ Knoxville, TN 37996-1300,
USA.}\\
\small{$\;^c\,$Department of Mathematics,
Hebei Normal University,}\\
\small{Shijiazhuang 050016, PRC}}
 \maketitle

\footnotetext[1]{Supported by an NSF grant.}
\footnotetext[2]{Supported by an NSA grant.}

\begin{abstract}
For a superprocess under a stochastic flow, we prove that it has a
density with respect to the Lebesgue measure for $d=1$ and is
singular for $d>1$. For $d=1$, a stochastic partial differential
equation is derived for the density. The regularity of the
solution is then proved by using Krylov's $L_p$-theory for linear
SPDE. A snake representation for this superprocess is established.
As applications of this representation, we prove the compact
support property for general $d$ and singularity of the process
when $d>1$.  \vspace*{.125in}


\vspace*{.125in}


\end{abstract}

\section{Introduction}
\setcounter{equation}{0}
\renewcommand{\theequation}{\thesection.\arabic{equation}}

Superprocesses under stochastic flows have been studied by many
authors since the work of Wang (\cite{W97},\cite{W1}) and
Skoulakis and Adler \cite{SA}. At an early stage, this problem was
studied as the high-density limit of a branching particle system
while the motion of each particle is governed by an independent
Brownian motion as well as by a common Brownian motion which
determines the stochastic flow. The limit is characterized by a
martingale problem whose uniqueness is established by a moment
duality. Before we go any further, let us introduce the model in
more detail.

 Let $b:\ \RR^d\to\RR^d$, $\si_1,\ \si_2:\
\RR^d\to\RR^{d\times d}$ be measurable functions. Let $W,\ B_1,\
B_2,\ \cdots$ be independent $d$-dimensional Brownian motions.
Consider a branching particle system performing independent binary
branching. Between branching times, the motion of the $i$th
particle is governed by the following stochastic differential
equation (SDE):
\begin{equation}\label{motion}
d\eta_i(t)=b(\eta_i(t))dt+\si_1(\eta_i(t))dW(t)+\si_2(\eta_i(t))dB_i(t).\end{equation}
It is proved by Skoulakis and Adler \cite{SA} that the
high-density limit $X_t$ is the unique solution to the following
martingale problem (MP): $X_0=\mu\in\cM_F(\RR^d)$, where $\cM_F(\RR^d)$
denotes the space of finite nonnegative measures on $\RR^d$ and for any
$\phi\in C^2_0(\RR^d)$,
\begin{equation}\label{eq1a3}
M_t(\phi)\equiv\left<X_t,\phi\right>-\left<\mu,\phi\right>-
\int^t_0\left<X_s,L\phi\right>ds
\end{equation}
is a continuous martingale with quadratic variation process
\begin{equation}\label{eq1a4}
\left<M(\phi)\right>_t=\int^t_0\left(\left<X_s,\phi^2\right>+
\left|\left<X_s,\si_1^T\nabla\phi\right>\right|^2\right)ds
\end{equation}
where
\[L\phi=\sum^d_{i=1}b^i\partial_i\phi
+\frac{1}{2}\sum^d_{i,j=1}a^{ij}\partial^2_{ij}\phi,\]
 $a^{ij}=\sum^d_{k=1}\sum^2_{\ell=1}\si^{ik}_{\ell}\si^{kj}_{\ell}$,
$\partial_i$ means the partial derivative with respect to the
$i$th component of $x\in\RR^d$, $\si^T_1$ is the transpose of the
matrix $\si_1$, $\nabla=(\partial_1,\cdots,\partial_d)^T$ is the
gradient operator and $\left<\mu,f\right>$ represents the integral
of the function $f$ with respect to the measure $\mu$. It was
conjectured in \cite{SA} that the conditional log-Laplace
transform of $X_t$ should be the unique solution to a nonlinear
stochastic partial differential equation (SPDE). Namely
\begin{equation}\label{eq1a2}
\EE_{\mu}\left(e^{-\left<X_t,f\right>}\Bigg| W\right)
=e^{-\left<\mu,y_{0,t}\right>}
\end{equation}
and
\begin{eqnarray}\label{eq1a1}
y_{s,t}(x)&=&f(x)+\int^t_s\left(L y_{r,t}(x)-y_{r,t}(x)^2\right)dr\nonumber\\
&&+\int^t_s\nabla^T y_{r,t}(x)\si_1(x)\hat{d}W(r)
\end{eqnarray}
where $\hat{d}W(r)$ represents the backward It\^o integral:
\[\int^t_sg(r)\hat{d}W(r)=\lim_{|\De|\to 0}\sum^n_{i=1}g\left(r_i\right)
\left(W\left(r_i\right)-W\left(r_{i-1}\right)\right)\] where
$\De=\{r_0,r_1,\cdots,r_n\}$ is a partition of $[s,t]$ and $|\De|$
is the maximum length of the subintervals.

This conjecture was confirmed by Xiong \cite{X1} under the
following conditions (BC) which will be assumed throughout this
paper: {\em $f\ge 0,\ b,\ \si_1,\ \si_2$ are bounded with bounded
first and second derivatives. $\si^T_2\si_2$ is uniformly positive
definite, $\si_1$ has third continuous bounded derivatives. $f$ is
of compact support.}

Making use of the conditional log-Laplace functional, the
long-term behavior of this process is studied in \cite{X2}. Also,
the model has been extended in that paper to allow infinite
measures $\mu\in\cM_{tem}(\RR^d)$, namely,
$\int_{\RR^d}e^{-\la|x|}\mu(dx)<\infty$ for some $\la>0$. We shall
assume $\mu\in\cM_{tem}(\RR^d)$ throughout this paper. A similar
model has been investigated by Wang \cite{W1} and Dawson et al
\cite{DLW} when the spatial dimension is 1. Further, in that case,
it is proved by Dawson et al \cite{DVW} that their process is
density-valued and solves a SPDE. The regularity of the solution
was left {\em open} in that article.

This paper is organized as follows: In Section 2, we establish a
snake representation for $X_t$. As immediate consequences to this
representation, we get the compact support property of $X_t$ (for
all $d$) and for $d>1$, $X_t$ takes values in the set of singular
measures. Then, for $d=1$, we prove in Section 3 that $X_t$ is
absolutely continuous with respect to Lebesgue measure and show
that the density $X(t,x)$ satisfies the following SPDE
\begin{equation}\label{main-eq}
\partial_t X=L^*X-\partial_x(\si_1X)\dot{W}_t
+\sqrt{X}\dot{B}_{tx} \end{equation} where $B$ is a Brownian sheet
and $L^*$ is the adjoint operator of $L$. The main result of this
paper is to show the H\"older continuity of $X(t,x)$.

Here is the main result.  First recall that
for $n\in\RR$ and $p\in[2,\infty)$,
$H^n_p$ is the space of Bessel potentials with norm
\[\|u\|_{n,p}=\|(I-\De)^{n/2}u\|_p.\]
\begin{thm}\label{main-thm}
Suppose that Condition (BC) is satisfied. Then\\
i) If $d>1$, then $X_t$ is singular a.s.\\
ii) If $d=1$, then $X_t$ is absolutely continuous with respect to
Lebesgue
measure and the density satisfies the SPDE (\ref{main-eq}).\\
iii) If in addition, $\mu$ satisfies $\mu\in
H^{\frac{1}{2}-\ep-2/p}_p$ with $\ep\in(0,\frac{1}{4})$ and
$p>\frac{1}{\ep}$ and also satisfies
\begin{equation}\label{cond4cont}
\sup_{t,x}\<\mu,\varphi_t(x-\cdot)\><\infty,
\end{equation}
 then the
density $X(t,x)$ is H\"older continuous in $x$ with index
$\frac{1}{2}-2\ep$ for (a.e.) $t$ a.s., where $\varphi_t(x)$ is
the density of a normal random variable with mean 0 and variance
$t$.
\end{thm}

Note that (\ref{cond4cont}) is satisfied if $\mu$ has bounded
density with respect to Lebesgue measure.

Suppose that we apply the usual integral equation as in \cite{wal86},
Chapter 3, for (\ref{main-eq}) in order to prove the H\"older continuity.
Then formally we have
\begin{eqnarray*}
X(t,x)&=&\int p_0(t,x,y)X(0,y)dy
+\int^t_0\int\si_1(y)X(s,y)\partial_yp_0(t-s,x,y)dy d W(s)\\
&&+\int^t_0\int\sqrt{X(s,y)}p_0(t-s,x,y)B(ds dy)
\end{eqnarray*}
where $p_0$ is the transition function of the Markov process with
generator $L$. However, the second term on the right hand side of
the above equation is about \[\int^t_0(t-s)^{-1/2}dW(s)\] which is
{\em not} convergent. Therefore, the convolution argument used by
Konno and Shiga \cite{KS} does not apply to our model. In Section
4, we freeze the nonlinear term in (\ref{main-eq}) and apply
Krylov's $L_p$-theory for linear SPDE to get the H\"older
continuity with index slightly less than $\frac{1}{2}$ for $X$.

Note that the SPDE in \cite{DVW} is (\ref{main-eq}) in current
paper with $\dot{W}_t$ replaced by a space-time
noise which is colored in space and white in time. The method of
this paper can be applied to that equation to prove the regularity
for its solution.

\section{Snake representation}
\setcounter{equation}{0}
\renewcommand{\theequation}{\thesection.\arabic{equation}}

In this section, we construct a path-valued process $\cY_t$ such
that the process $X_t$ can be represented according to this
process. Then, as an easy application of this representation, we
derive the properties for $X_t$.

For the convenience of the reader, we recall some basic
definitions and facts taken from Le Gall \cite{LeG}. Let $\ze\ge
0$ and let $f$ be a continuous function from $\RR_+$ to $\RR^d$
such that $f(s)=f(\ze)$, $\forall\ s\ge \ze$. We call such pair
$(f,\ze)$ a stopped path with $\ze$ being the lifetime of the
path. We denote the collection of all stopped paths by $\WW$. For
$(f,\ze),\;(f',\ze')\in\WW$, define a distance
\[\de((f,\ze),(f',\ze'))=\sup_{s\ge 0}|f(s)-f'(s)|+|\ze-\ze'|.\]
Then $(\WW,\de)$ is a Polish space. In \cite{LeG}, Le Gall
constructed a continuous time-homogeneous strong Markov process
$(\cZ_t,\ze_t)$ taking values on $\WW$. $\ze_t$ is a
one-dimensional reflecting Brownian motion. Given $\ze$, the
process $\cZ$ has the following property:  for all $r<t$, and for
all $s\le m_{r,t}:=\inf_{r\leq u\leq t}\ze_u$ we have
$\cZ_r(s)=\cZ_t(s)$.  Furthermore, given $m_{r,t}$ and
$\cZ_r(m_{r,t})$, the processes $\cZ_r(s): s\geq m_{r,t}$ and
$\cZ_t(s): s\geq m_{r,t}$ are conditionally independent Brownian
motions with lifetimes $\ze_r$ and $\ze_t$ respectively.

Denote the strong solution to the SDE
\[d\eta(t)=b(\eta(t))dt+\si_1(\eta(t))dW(t)+\si_2(\eta(t))dB(t)\]
by $\eta(t)=F(t,W,B)$. Define the following path-valued process
\[\cY_t(s)=F(s,W,\cZ_t)\]
with the life-time process $\ze_t$.

\begin{lemma}
$(\cY_t,\ze_t)$ is a continuous $\WW$-valued process.
\end{lemma}
Proof: Note that for all $r<t$ and for all $s<m_{r,t}$, we have
$\cY_r(s)=\cY_t(s)$.  Furthermore, for
given $\cY_r(m_{r,t})$, the processes $\cY_r(s): s\geq m_{r,t}$
and $\cY_t(s): s\geq m_{r,t}$ are the motions of two particles
(say, $\eta_1$ and $\eta_2$) given as in the introduction with
lifetimes $\ze_r$ and
$\ze_t$ starting from the same position $\cY_r(m_{r,t})$. A simple
application of Burkholder's inequality gives
\[\EE\left[\sup_{m\le s\le M}|\eta_1(s)-\eta_2(s)|^k\right]\le K|M-m|^{k/2},\]
where $m=m_{r,t}$ and $M=\ze_r\vee\ze_t$. Denote by $\EE^{\ze}$
the conditional expectation given $\ze$. Then
\begin{eqnarray*}
\EE\left[\sup_{s\ge 0}|\cY_r(s)-\cY_t(s)|^k\right]
&=&\EE\left[\EE^{\ze}\left\{\sup_{s\ge m_{r,t}}|\cY_r(s)-\cY_t(s)|^k\right\}\right]\\
&\le&\EE \left[K|\ze_r+\ze_t-2m_{r,t}|^{k/2}\right]\\
&\le&K\EE\left[\sup_{s\in[r,t]}|\ze_s-\ze_r|^{k/2}\right]\\
&\le&K|t-r|^{k/4}.
\end{eqnarray*}
The conclusion follows from Kolmogorov's criteria by taking $k>4$;
see \cite{wal86} for Kolmogorov's criteria.
\qed

\begin{thm}
\begin{equation}\label{snake}
X_t(f)=\int^{\tau}_tf(\cY_s(\ze_s))d\ell^t_s\end{equation} where
$\ell^t$ is the local time process of $\ze$ at level $t$ and
\[\tau=\inf\{s:\;\ell^0_s\ge 1\}.\]
\end{thm}
Proof: Fix a parameter $h>0$. For every $t\ge 0$, denote by
$[a^1_t,b^1_t]$, $[a^2_t,b^2_t],\;\cdots,\; [a^{N_t}_t,b^{N_t}_t]$
the excursion intervals of $(\ze_s)_{0\le s\le\tau}$ above level
$t$, corresponding to excursions of height greater than $h$. Set
\[X^h_t=2h\sum^{N_t}_{i=1}\de_{\cY_{a^i_t}(t)}.\]
Then $X^h_t$ is the measure-valued process corresponding to the
branching particle system described as follows: At time $t=0$, we
have $N_0$ particles in $\RR^d$ with Poisson random measure with
intensity measure $h^{-1}\mu$. The particles then move according
to (\ref{motion}) with common $W$ and independent $B_i$'s. Each of
them has a finite lifetime (independent of others) which is
exponential with mean $h$. When a particle dies, it gives rise to
either 0 or 2 new particles with probability $\frac{1}{2}$. The
new particles start from the position of the their father.
As in the proof of Theorem 2.1 in \cite{LeG}, by the well-known
approximation of Brownian local time by upcrossing numbers, we
have that $X^h_t$ converges weakly to $X_t$, where $X_t$ is given
by the right hand side
of (\ref{snake}). \qed

As an application of the snake representation, we have the
following immediate consequence.

\begin{cor}
If $\mu$ is a finite measure, then for any $t>0$, $X_t$ has
compact support a.s.
\end{cor}
Proof: By the snake representation, there exists a finite set $I$
such that
\[\<X_t,f\>=\sum_{i\in I}\int^{\tau_i}_0f(\hat{\cY}^i_s)
d\ell^t_s(\ze^i)\] where $\hat{\cY}^i_s$ is the tip of the $i$th
snake. It is not hard to show that $\hat{\cY}^i_s$ is continuous and
hence, for any $t_0>0$,
\begin{equation}\label{supp}
\bigcup_{t\ge t_0}\mbox{supp}(X_t)\subset\overline{\bigcup_{i\in
I}{\rm Range}\(\hat{\cY}^i\)}=\bigcup_{i\in I}\{\hat{\cY}^i_s:\;0\le
s\le\tau_i\}
\end{equation}
is compact. \qed

To consider the case for $\mu$ being $\si$-finite, the following
conditional martingale problem (CMP) is useful. The following
lemma was proved in \cite{X2}.
\begin{lemma}
i) If $X_t$ is the solution to MP, then there exists a Brownian motion $W_t$
such that for any $\phi\in C^2_0(\RR^d)$,
\begin{equation}\label{cmp1}
N_t(\phi)\equiv\left<X_t,\phi\right>-\left<\mu,\phi\right>-
\int^t_0\left<X_s,L\phi\right>ds-\int^t_0\left<X_s,\si_1^T\nabla\phi\right>dW_s
\end{equation}
is a continuous $(\PP,\cG_t)$-martingale with quadratic variation
process
\begin{equation}\label{cmp2}
\left<N(\phi)\right>_t=\int^t_0\left<X_s,\phi^2\right>ds
\end{equation}
where $\cG_t=\cF_t\vee\cF^W_{\infty}$.

ii) If $X_t$ is a solution to CMP, then it is a solution to MP.
\end{lemma}

As another application of the snake representation, we have
\begin{cor}\label{cor-2}
If $d\ge 2$, then $X_t$ is singular.
\end{cor}
Proof: If $\mu$ is finite and $d>1$, it follows from (\ref{supp})
the support is of Lebesgue measure $0$ since
$\{\hat{\cY}^i_s:\;0\le s\le\tau_i\}$ is a continuous
(one-dimensional) curve in $\RR^d$. If $\mu$ is
$\si$-finite, we can take $\mu=\sum^{\infty}_{n=1}\mu^n$ with
$\mu^n$ finite. Construct the solution $X^n_t$ to CMP with the
same $W$ and with initial $\mu^n$, $n=1,2,\cdots$. Then
\[X_t=\sum^{\infty}_{n=1}X^n_t\]
is the solution to CMP with initial $\mu$. Then
$\mbox{supp}(X^n_t)$ has Lebesgue measure 0 and hence, so does the
support of $X_t$. This implies that $X_t$ is a singular measure
a.s. \qed

\section{SPDE for $d=1$}
\setcounter{equation}{0}
\renewcommand{\theequation}{\thesection.\arabic{equation}}

In this section, we prove that $X_t$ has a density which satisfies
the SPDE (\ref{main-eq}) whose mild form is
\begin{eqnarray}\label{spde}
\left<X_t,f\right>&=&\left<\mu,f\right>+\int^t_0\left<X_s,Lf\right>ds
+\int^t_0\left<X_s,\si_1 f'\right>dW_r\nonumber\\
&&+\int^t_0\int_{\RR}\sqrt{X_s(x)}f(x)B(ds dx).
\end{eqnarray}

Let $p_0(t,x,y)$ and $q_0(t,(x_1,x_2),(y_1,y_2))$ be the
transition density functions of the Markov processes $\eta_1(t)$
and $(\eta_1(t),\eta_2(t))$ respectively. By Theorem 1.5 of
\cite{X1}, we have
\begin{equation}\label{1stmom}
\EE\Big[\<X_t,f\>\Big]=\int_{\RR^2}f(y)p_0(t,x,y)dy\mu(dx)
\end{equation}
 and
\begin{eqnarray}\label{2ndmom}
\lefteqn{\EE\Big[\<X_t,f\>\<X_t,g\>\Big]  }\\
&=&\int_{\RR^4}f(y_1)g(y_2)q_0(t,(x_1,x_2),(y_1,y_2))dy_1dy_2\mu(dx_1)\mu(dx_2)\nonumber\\
&&+2\int^t_0ds\int_{\RR^4}p_0(t-s,z,y)f(z_1)g(z_2)q_0(s,(y,y),(z_1,z_2))dz_1dz_2dy\mu(dz).
\nonumber
\end{eqnarray}

\begin{thm}
If $\mu(\RR)<\infty$, then $X_t\in H_0\equiv L^2(\RR)$ a.s.
\end{thm}
Proof: Take $f=p_0(\ep,x,\cdot)$ and $g=p_0(\ep',x,\cdot)$ in
(\ref{2ndmom}). Note that as $\ep,\ep'\to 0$,
\begin{eqnarray*}
&&\int_{\RR^2}p_0(\ep,x,z_1)p_0(\ep',x,z_2)p_0(t-s,z,y)q_0(t,(y,y),(z_1,z_2))dz_1dz_2\\
&\to& p_0(t-s,z,y)q_0(t,(y,y),(x,x)).\end{eqnarray*} Note that by
Theorem 6.4.5 in Friedman \cite{F}, we have
\[p_0(\ep,x,y)\le c\varphi_{c'\ep}(x-y),\]
\[q_0(s,(y,y),(z_1,z_2))\le c\varphi_{c's}(y-z_1)\varphi_{c's}(y-z_2)\]
where $\varphi_t(x)$ is the normal density with mean 0 and
variance $t$ (introduced earlier). Note that
$c'$ is a constant which is usually greater than 1. Since it does
not play an essential role, to simplify the
notations, we assume $c'=1$ throughput the rest
of this paper. Hence,
\begin{eqnarray*}
\lefteqn{\int_{\RR^2}p_0(\ep,x,z_1)p_0(\ep',x,z_2)p_0(t-s,z,y)q_0(s,(y,y),(z_1,z_2))dz_1dz_2}\\
&\le&c\int_{\RR^2}\varphi_{\ep}(x-z_1)\varphi_{\ep'}(x-z_2)\varphi_{t-s}(z-y)\varphi_s(y-z_1)\varphi_s(y-z_2)dz_1dz_2\\
&=&c\varphi_{s+\ep}(x-y)\varphi_{s+\ep'}(x-y)\varphi_{t-s}(z-y).
\end{eqnarray*}
As
\begin{eqnarray*}
\lefteqn{\lim_{\epsilon,\epsilon'\to0}
\int^T_0dt\int dx\int^t_0ds\int_{\RR^2}\varphi_{s+\ep}(x-y)\varphi_{s+\ep'}(x-y)\varphi_{t-s}(z-y)dy
\mu(dz)}\\
&=&\lim_{\epsilon,\epsilon'\to0}\int^T_0dt\int^t_0ds\varphi_{2s+\ep+\ep'}(0)\mu(\RR)\\
&=&\int^T_0dt\int^t_0ds\varphi_{2s}(0)\mu(\RR)      \\
&=&\int^T_0dt\int dx\int^t_0ds\int_{\RR^2}\varphi_{t-s}(z-y)\varphi_{s}(x-y)\varphi_{s}(x-y)dy\mu(dz),
\end{eqnarray*}
by the dominated convergence theorem, we see that as $\ep,\ep'\to 0$,
\begin{eqnarray*}
&&\int^T_0dt\int dx\int^t_0ds\int_{\RR^4}p_0(t-s,z,y)p_0(\ep,x,z_1)
p_0(\ep',x,z_2)q_0(s,(y,y),(z_1,z_2))dz_1dz_2dy\mu(dz)\\
&\to&\int^T_0dt\int
dx\int^t_0ds\int_{\RR^2}p_0(t-s,z,y)q_0(t,(y,y),(x,x))dy\mu(dz).
\end{eqnarray*}
Similarly, we have
\begin{eqnarray*}
&&\int^T_0dt\int dx\int_{\RR^4}p_0(\ep,x,y_1)p_0(\ep',x,y_2)
q_0(t,(x_1,x_2),(y_1,y_2))dy_1dy_2\mu(dx_1)\mu(dx_2)\\
&\to&\int^T_0dt\int
dx\int_{\RR^2}q_0(t,(x_1,x_2),(x,x))\mu(dx_1)\mu(dx_2).
\end{eqnarray*}
Hence
\begin{eqnarray*}
&&\int^T_0dt\int
dx\EE\(\<X_t,p(\ep,x,\cdot\>\<X_t,p(\ep',x,\cdot)\>\)\\
 &\to&\int^T_0dt\int dx
\int_{\RR^2}q_0(t,(x_1,x_2),(x,x))\mu(dx_1)\mu(dx_2)\\
&&+\int^T_0dt\int
dx\int^t_0ds\int_{\RR^2}p_0(t-s,x,y)q_0(t,(y,y),(x,x))dy\mu(dx).
\end{eqnarray*}
From this, we can show that $\{\<X_t,p_0(\ep,x,\cdot)\>:\;\ep>0\}$
is a Cauchy sequence in $L^2(\Om\times[0,T]\times \RR)$. This
implies the existence of the density $X_t(x)$ of $X_t$ in
$L^2(\Om\times[0,T]\times \RR)$. \qed

Next theorem considers infinite measure.

\begin{thm}
If $\mu\in\cM_{tem}(\RR^d)$, then $X_t$ has a density $X_t(x)$.
\end{thm}
Proof: If $\mu$ is $\si$-finite, we can construct $X^n$ with
$X^n_0=\mu^n$ being finite as those in the proof of Corollary
\ref{cor-2}. Then
\[X_t=\sum^{\infty}_{n=1}X^n_t\]
is the solution to CMP with initial $\mu$. Let
\begin{equation}\label{density}
X_t(x)=\sum^{\infty}_{n=1}X^n_t(x). \end{equation}
 By
(\ref{1stmom}), we have
\[\EE X^n_t(x)=\int_{\RR}p_0(t,y,x)\mu^n(dy).\]
As \[p_0(t,x,y)\le c\varphi_t(x-y)\le c(t,\la,x)e^{-\la|y|},\] for
any $\la>0$, we have
\begin{eqnarray*} \EE
\sum^{\infty}_{n=1}X^n_t(x)&=&\sum^{\infty}_{n=1}\int_{\RR}p_0(t,y,x)\mu^n(dy)\\
&=&\int_{\RR}p_0(t,y,x)\mu(dy)<\infty.\end{eqnarray*} Hence,
$X_t(x)$ is well-defined by (\ref{density}). It is then easy to
show that $X_t(x)dx=X_t(dx)$. \qed

Finally, we derive the SPDE satisfied by the density.

\begin{thm}
If $d=1$, then $X_t$ is the (weak) unique solution to the SPDE
(\ref{spde}).
\end{thm}
Proof: Note that $N_t(\phi)$ is a continuous
$(\PP,\cG_t)$-martingale with quadratic variation process
\[\<N(\phi)\>_t=\int^t_0\int_{\RR}\(\sqrt{X_s(x)}\phi(x)\)^2dxds.\]
By the martingale representation theorem (\cite{KX}, Theorem
3.3.5), there exists an $L^2(\RR)$-cylindrical Brownian motion
$\tilde{B}$ on an extension of $(\Om,\cF,\cG_t,\PP)$ such that
\[N_t=\int^t_0\<\sqrt{X_s},d\tilde{B}_s\>_{L^2(\RR)}.\]
There exists a standard Brownian sheet $B$  such that
 \[\tilde{B}_t(h)=\int^t_0\int_{\RR} h(x)B(ds dx),\qquad\forall h\in
 L^2(\RR).\] Therefore,
\[N_t(\phi)=\int^t_0\int_{\RR}\sqrt{X_s(x)}\phi(x)B(ds dx).\]
As $B$ is a Brownian sheet on an extension of $\cG_t$, it is easy
to show that $B$ is independent of $W$. \qed

\section{H\"older Continuity}
\setcounter{equation}{0}
\renewcommand{\theequation}{\thesection.\arabic{equation}}

This section is devoted to the proof of the main result: Theorem
\ref{main-thm} (iii). Namely, in this section, we consider the
regularity of the solution to the nonlinear SPDE (\ref{main-eq}).
We use the linearization and Krylov's $L_p$-theory for linear
SPDE.

We will paraphrase the condition (BC) to find some reasonable
assumptions for $\sigma_1,\sigma_2,b $ to make our regularity
argument easy. Note that these functions are scalar functions
since we are dealing with the situation $d=1$. Therefore, we have
$L=\frac{1}{2}a\partial_{xx}+b\partial_x$ and
$L^*=\frac{1}{2}a\partial_{xx}+(a'-b)\partial_x+(\frac{1}{2}a''-b')$.

We start by defining some basic spaces. We denote
\[
[f]_0=\sup_{x\in \RR}|f(x)|,\quad
[f]_{\gamma}=\sup_{x\ne y}\frac{|f(x)-f(y)|}{|x-y|^{\gamma}}
\]
for $\gamma\in (0,1]$. Using this notation, we define
\begin{eqnarray*}
\|f\|_{C^{0,\gamma}}=[f]_0+[f]_{\gamma},\quad
\|f\|_{C^{1,\gamma}}=[f]_0+[f']_0+[f']_{\gamma}\\
\|f\|_{C^{1}}=[f]_0+[f']_0,\quad
\|f\|_{C^2}=[f]_0+[f']_0+[f'']_0
\end{eqnarray*}
assuming that $f'$ or $f''$ exist if they appear in the corresponding
definition.  Then we define the Banach spaces :
\begin{eqnarray*}
C^{0,\gamma}=\{f\;:\;\|f\|_{C^{0,\gamma}}<\infty\},\quad
C^{1,\gamma}=\{f\;:\;\|f\|_{C^{1,\gamma}}<\infty\}\\
C^{1}=\{f\;:\;\|f\|_{C^1}<\infty\},\quad
C^{2}=\{f\;:\;\|f\|_{C^2}<\infty\}.
\end{eqnarray*}
\begin{rem}\label{Zygmund}
Zygmund spaces $C^{0,\gamma},C^{1,\gamma}$ are the usual H\"older
spaces if $\gamma\in (0,1)$. It is easy to see that we have
$\|f\|_{C^{0,\gamma}}\le 2\|f\|_{C^{0,1}},\|f\|_{C^{1,\gamma}}\le
2\|f\|_{C^{1,1}}$ and $\|f\|_{C^{0,1}}\le
\|f\|_{C^{1}},\|f\|_{C^{1,1}}\le \|f\|_{C^{2}}$ when $f'$ or $f''$
exists.
\end{rem}

Now, we state assumptions on $\sigma_1,\sigma_2,b$. First, our
condition (BC) gives us the following assumption :
\begin{equation}\label{sigmas}
\sigma_1,\sigma_2\in C^{2},\quad b\in C^{1}
\end{equation}
which, in particular, implies $a=\sigma^2_1+\sigma^2_2\in C^{2}$.
We also assume that
\begin{equation}\label{bound of norms}
\delta\le\frac{1}{2}a,\frac{1}{2}\sigma^2_2\le K,\quad
\|\sigma_1\|_{C^{2}},\|\sigma_2\|_{C^{2}},\|b\|_{C^{1}}\le K
\end{equation}
for some positive constants $\delta,K$.

Next, we recall the basic definitions of some function spaces
defined  in \cite{Krylov1999}. In addition to the definition about
space of Bessel potentials in the Theorem \ref{main-thm}, we also
define the following : for $n\in\RR$ and $p\in[2,\infty)$ let
$H^n_p(l_2) $ be the space with norm
\[
\|g\|_{n,p}=\Big\| |(I-\De)^{n/2}g|_{l_2}\Big\|_p
\]
for $l_2 $-valued functions $g=\{g^k\} $. Then we define
\[
\HH^n_p(T)=L_p(\Om\times[0,T],\cP,H^n_p)\quad
\HH^n_p(T,l_2)=L_p(\Om\times[0,T],\cP,H^n_p(l_2))
\]
where $\cP$ is the predictable $\si$-field. We denote
$\bb{L}_p(T)= \HH^0_p(T)$. Let $\{w^k_t:k=1,2,\ldots\} $ be a
family of independent one-dimensional Brownian motions.

We say $u\in\cH^n_p(T)$ if $\partial_{xx}u\in \HH^{n-2}_p(T)$ and
$u(0,\cdot)\in L_p(\Om,H^{n-2/p}_p)$ and  there exists $(f,g)\in
\HH^{n-2}_p(T)\times\HH^{n-1}_p(T,l_2)$ such that $\forall\phi\in
C^{\infty}_0(\RR)$, (a.s.)
\[\<u_t,\phi\>=\<u_0,\phi\>+\int^t_0\<f_s,\phi\>ds+\sum^{\infty}_{k=0}\int^t_0\<g^k_s,\phi\>dw^k_s\]
holds for all $t\le T $. We denote
\[\|u\|_{\cH^n_p(T)}=\|\partial_{xx}u\|_{\HH^{n-2}_p(T)}
+\|f\|_{\HH^{n-2}_p(T)}+\|g\|_{\HH^{n-1}_p(T,l_2)}+\(\EE\|u_0\|^p_{n-2/p,p}\)^{1/p}\]
Reader can find motivation of this definition and detailed remarks
in \cite{Krylov1999}.

Now, we fix $\ep\in (0,\frac{1}{4})$ and proceed to the
 {\em Proof of Theorem \ref{main-thm} (iii)} :
First, we freeze the nonlinear term of SPDE (\ref{main-eq}) and
consider the following auxiliary linear SPDE for $Y_t(x)$:
\begin{equation}\label{lspde}
\left\{\begin{array}{l}\partial_tY=L^*Y+\sqrt{X}\dot{B}_{tx}\\
Y_0=\mu
\end{array}\right.
\end{equation}
where we assume $\mu\in H^{\frac{1}{2}-\ep-2/p}_p$.

 Then $Z=X-Y$  satisfies
\begin{equation}\label{spde4z}
\left\{\begin{array}{l}\partial_tZ=L^*Z-\(\partial_x(\si_1Z)+\partial_x(\si_1Y)\)\dot{W}_t\\
Z_0=0.
\end{array}\right.
\end{equation}

We apply Theorem 8.5 of \cite{Krylov1999} to (\ref{lspde}). To do
this we need the coefficients of $L^*$ and $\sqrt{X}$ to satisfy
\[
\|a\|_{C^{1,1}}<\infty,\quad\|a'-b\|_{C^{0,1}}<\infty,
\quad[\frac{1}{2}a''-b']_0<\infty,\quad\|\sqrt{X}\|_{\bb{L}_p(T)}<\infty.
\]
In fact, we have
\[
\|a\|_{C^{1,1}}\le K,\quad\|a'-b\|_{C^{0,1}}\le 2K,\quad[\frac{1}{2}a''-b']_0\le 2K
\]
by our assumptions (\ref{sigmas}) and (\ref{bound of norms}) and
Remark \ref{Zygmund}. We will prove
$\|\sqrt{X}\|_{\bb{L}_p(T)}<\infty$ later and take this for
granted in this proof.

Now, by Theorem 8.5 of \cite{Krylov1999} to (\ref{lspde}) and the
first assertion of Lemma 8.4 and the fact that $\mu$ is nonrandom,
we have a unique solution $Y$ in $\cH^{\frac{1}{2}-\ep}_p(T)$ with
estimate
\begin{eqnarray}\label{bound4y}
\|Y\|_{\cH^{\frac{1}{2}-\ep}_p(T)}
&\le&N(\|\sqrt{X}\|_{\bb{L}_p(T)}+\|\mu\|_{\frac{1}{2}-\ep-2/p,p})
\end{eqnarray}
where $N$ depends only on $\ep,p,\delta,K,T$.

Now we use Theorem 5.1 in \cite{Krylov1999} for equation
(\ref{spde4z}) above with $n=-\frac{3}{2}-\ep\in
(-2,-\frac{3}{2})$. Note $\partial_x(\si_1Z)=\si_1\partial_x
Z+\partial_x\si_1 Z$ . If we read \cite{Krylov1999} carefully, we
can see that the following conditions are required :
\begin{enumerate}
\item[(i)] \[
\delta\le
\frac{1}{2}a-\frac{1}{2}\sigma^2_1(=\frac{1}{2}\sigma^2_2)\le K_1
\]
for some positive $\delta,K_1 $.
\item[(ii)] $a,\sigma_1 $ are Lipschitz continuous with Lipschitz constant $K_1 $.
\item[(iii)] $a\in C^{1,\gamma_1},\sigma_1\in C^{0,\gamma_2} $  for some
$\gamma_1,\gamma_2\in (0,1)$ and
$\|a\|_{C^{1,\gamma_1}}+\|\sigma\|_{C^{0,\gamma_2}}\le K_1 $
\item[(iv)] $\|a'-b\|_{C^{0,\gamma_3}}+[\frac{1}{2}(a''-b')]_0+[\partial_x\sigma_1]_0\le
K_1 $ for some $\gamma_3\in (0,1)$.
\item[(v)] $\partial_x(\sigma_1Y)\in\bb{H}^{n+1}_p(T) \;(=\bb{H}^{-\frac{1}{2}-\ep}_p(T)$).
\end{enumerate}
  But, conditions (i) through (iv) are satisfied under
(\ref{sigmas}) and (\ref{bound of norms}) and Remark
\ref{Zygmund}. Note that we can take some constant multiple of
$K^2 $ as $K_1 $. On the other hand, (v) is also satisfied. For
\begin{eqnarray}\label{bound4z}
\left\|\partial_x(\si_1Y)\right\|_{\HH^{-\frac{1}{2}-\ep}_p(T)}
&\le& N\|\si_1Y\|_{\HH^{\frac{1}{2}-\ep}_p(T)}\label{operator}\\
 &\le&
N\|\si\|_{C^{0,{\frac{1}{2}-\ep+\frac{1}{4}}}}\|Y\|_{\HH^{\frac{1}{2}-\ep}_p(T)}\label{co}\\
 &\le&
N\|\si\|_{C^1}\|Y\|_{\HH^{\frac{1}{2}-\ep}_p(T)}\label{free}\\
 &\le&
N\|Y\|_{\cH^{\frac{1}{2}-\ep}_p(T)}\label{cH}\\
&\le&
N\|\sqrt{X}\|_{\bb{L}_p(T)}+N\|\mu\|_{\frac{1}{2}-\ep-2/p,p}<\infty
.\label{last}
\end{eqnarray}
(\ref{operator}) follows the observation
$\partial_x=\partial_x(I-\Delta)^{-1/2}(I-\Delta)^{1/2} $ and the
boundness of the operator $\partial_x(I-\Delta)^{-1/2}$.
(\ref{co}) follows Lemma 5.1 (i) in \cite{Krylov1999}. Up to this
step, $N $ only depends on $\ep,p$. Note that
$\frac{1}{2}-\ep+\frac{1}{4} $ is still in $(0,1) $ since
$\frac{1}{2}-\ep\in (\frac{1}{4},\frac{1}{2} )$. Hence, we have
(\ref{free}) by (\ref{bound of norms}) and Remark \ref{Zygmund}.
(\ref{cH}) follows Theorem 3.7 in \cite{Krylov1999} and $N $
depends only on $\ep,p,K,T $ now. Finally, (\ref{bound4y}) gives
us (\ref{last}) with $N=N(\ep,p,\delta,K,T ) $.

Therefore, we have a unique solution $Z $ in $
\cH^{\frac{1}{2}-\ep}_p(T)$ with
\begin{equation}\label{bound4z2}
\|Z\|_{\cH^{\frac{1}{2}-\ep}_p(T)}\le
N\left\|\partial_x(\si_1Y)\right\|_{\HH^{-\frac{1}{2}-\ep}_p(T)}
\le N\|\sqrt{X}\|_{\bb{L}_p(T)}+N\|\mu\|_{\frac{1}{2}-\ep-2/p,p}
\end{equation}
where $N=N(\ep,p,\delta,K,T) $.

Thus, combining (\ref{bound4y}) and (\ref{bound4z2}), we have
$X=Y+Z \in \cH^{\frac{1}{2}-\ep}_p(T)$ with estimate
\begin{equation}\label{bound2x}
\|X\|_{\cH^{\frac{1}{2}-\ep}_p(T)}\le
N\|\sqrt{X}\|_{\bb{L}_p(T)}+N\|\mu\|_{\frac{1}{2}-\ep-2/p,p}.
\end{equation}
By the embedding Theorem 7.1 in \cite{Krylov1999}, this implies
\[\left(E\int^T_0\|X_t\|^p_{C^{{\frac{1}{2}}-\ep-\frac{1}{p}}}dt\right)^{1/p}\le N\|X\|_{\cH^{\frac{1}{2}-\ep}_p(T)}\le
N\|\sqrt{X}\|_{\bb{L}_p(T)}+N\|\mu\|_{\frac{1}{2}-\ep-2/p,p}.\]
 So, for large $p>\frac{1}{\ep}$
, we have
\[\|X_t\|_{C^{{\frac{1}{2}}-2\ep}}<\infty\]
for (a.e.) $t\in [0,T]$ a.s.. we are done with the proof. \qed

Finally, we use the moment dual to prove that
\begin{equation}\label{nth-mom}
\EE\int^T_0\int_{\RR}X(t,x)^ndx dt<\infty
\end{equation}
for all $n\in\NN$.

Let $n_t$ be a pure-death Markov chain with $n_0=0$ and, at a rate
$\frac{1}{2}n(n-1)$, jumps from $n$ to $n-1$. Let
$0=\tau_0<\tau_1<\cdots<\tau_{n-1}$ be the jump times. Let
$f_0=\de_x^{\otimes n}$ and for $t<\tau_1$,
$f_t(y)=p^n_0(t,(x,\cdots,x),y),\;\forall y\in\RR^n$ where $p^n_0$
is the transition function of the $n$-dimensional diffusion
$(\eta_1(t),\cdots,\eta_n(t))$. For $f\in C(\RR^n)$, let
$G_{ij}f\in C(\RR^{n-1})$ be given by
\[G_{ij}f(y_1,\cdots,y_{n-2},y_{n-1})=f(y_1,\cdots,y_{n-1},\cdots,y_{n-1},\cdots,y_{n-2})\]
where $y_{n-1}$ is at $i$th and $j$th position. Let
\[f_{\tau_1}=\Ga_1f_{\tau_1-}\]
where $\Ga_1$ is a random element taking values in $\{G_{ij}:1\le
i<j\le n\}$ uniformly. We continue this procedure to get the
process $f_t$. Replace $f_0$ by a smooth function $f_0^k\ge 0$
approximating $f_0$. Denote the process constructed above with
$f_0^k$ in place of $f_0$ by $f_t^k$. Similar to Theorem 11 in
Xiong and Zhou \cite{XZ}, we have
\[\EE\<X_t^{\otimes
n},f_0^k\>=\EE\(\<\mu^{\otimes
n_t},f^k_t\>\exp\(\frac{1}{2}\int^t_0n_s(n_s-1)ds\)\).\] Taking
limits and using Fatou's lemma, we have
\begin{eqnarray*}
\EE X(t,x)^n
&\le&\EE\(\<\mu^{\otimes n_t},f_t\>\exp\(\frac{1}{2}\int^t_0n_s(n_s-1)ds\)\)\\
&\le&\exp\(\frac{1}{2}n(n-1)t\)\sum^n_{i=1}\EE\(\<\mu^{\otimes
n_t},f_t\>1_{\tau_{i-1}\le t<\tau_i}\).
\end{eqnarray*}
Let $i=3$. Then
\begin{eqnarray*}
f_t(x_1,\cdots,x_{n-2})
&\le&c\int_{\RR^{n-2}}\Pi^{n-2}_{i=1}\varphi_{t-\tau_2}(x_i-y_i)\Ga_2f_{\tau_2-}(y)dy\\
&\le&c\int_{\RR^{n-2}}\Pi^{n-2}_{i=1}\varphi_{t-\tau_2}(x_i-y_i)
\sum_{1\le k<\ell}\frac{2}{(n-2)(n-3)}\\
&&\qquad\qquad f_{\tau_2-}(y_1,\cdots,y_{n-2},\cdots,y_{n-2},\cdots,y_{n-3})dy\\
&\le&c\int_{\RR^{n-2}}\Pi^{n-2}_{i=1}\varphi_{t-\tau_2}(x_i-y_i)
\sum_{1\le k<\ell}\frac{2}{(n-2)(n-3)}\\
&&\qquad\qquad
\int_{\RR^{n-1}}\Pi^{n-1}_{j=1}\varphi_{\tau_2-\tau_1}(y_j-z_j)\varphi_{\tau_1}(z_1-x)
\cdots\\
&&\qquad\qquad\cdots\varphi_{\tau_1}(z_{n-2}-x)\varphi_{\tau_1}(z_{n-1}-x)^2dz.
\end{eqnarray*}
Thus
\begin{eqnarray*}
\<\mu^{\otimes n-2},f_t\>
&\le&c\int_{\RR}\varphi_{t-\tau_2}(x_{n-2}-y_{n-2})\mu(dx_{n-2})\int_{\RR}dy_{n-2}
\sum_{1\le k<\ell}\frac{2}{(n-2)(n-3)}\\
&&\qquad\qquad\varphi_{\tau_2-\tau_1}(y_{n-2}-z_k)\varphi_{\tau_2-\tau_1}(y_{n-2}-z_{\ell})\\
&&\qquad\qquad \int_{\RR^{n-1}}\varphi_{\tau_1}(z_1-x)\cdots
\varphi_{\tau_1}(z_{n-2}-x)\varphi_{\tau_1}(z_{n-1}-x)^2\\
&\le&c\int_{\RR}\varphi_{t-\tau_2}(x_{n-2}-y_{n-2})\mu(dx_{n-2})\int_{\RR}dy_{n-2}
\frac{1}{\sqrt{\tau_1(\tau_2-\tau_1)}}\varphi_{\tau_2}(y_{n-2}-x).
\end{eqnarray*}
Therefore
\[\int_{\RR}\EE\<\mu^{\otimes n_t},f_t\>1_{\tau_2\le
t<\tau_3}dx \le
c\mu(\RR)\EE\frac{1}{\sqrt{\tau_1(\tau_2-\tau_1)}}<\infty.\] The
other terms can be proved similarly. This finishes the proof of
Theorem \ref{main-thm}.
\\[1cm]
 {\bf Acknowledgement:} Most of this work was done during the
second author's visit of University of Tennessee and the third
author's visit of University of Rochester. Financial support and
hospitality from both institutes are appreciated.

\end{document}